# A Multi-Period Water Network Planning for Industrial Parks; Impact of Design Periods on Park's Flexibility

Saman Taheri [1]

[1] Department of Mechanical Engineering, Purdue School of Engineering and Technology, Indiana University-Purdue University Indianapolis, IN 46202, USA

*staheri@iu.edu

**Abstract:** The central goal of this study is to advance the state-of-the-art in the field of water network management for industrial ecologies. Considering the global water shortage and environmental regulations, the reduction of water consumption has become a necessity within the industrial sector. To find the optimal water management strategies in an eco-industrial park (EIP), a refined superstructure is proposed. An MINLP mathematical model is then developed, based on the superstructure, to synthesize water and wastewater streams within the network, including nearby companies. Surveying two case studies, each studied in two different scenarios, showed that the cost objective functions deteriorates with increase in number of design periods. Results also provide some evidence of a trade-off between water network flexibility and performance, and vice versa.

**Nomenclature**

*Sets and Indices*

| | |
|---|---|
| $c$ | Water contaminants |
| $i$ | Process sources including $i \in \{1,2,...,n\}$ |
| $j$ | Process sinks including $j \in \{1,2,...,m\}$ |
| $s$ | Regenerators including $s \in \{1,2,...,S\}$ |
| $t$ | Periods including $t \in \{1,2,...,r\}$ |

*Variables*

| | |
|---|---|
| $C_{exp}$ | The contaminant level of exported flowrates |
| $C_{inv}$ | The annualized investment cost of the park |
| $C_{mix-CH}$ | The contaminant level of the CH flowrate |
| $C_{op}$ | The annaulized operation cost of the park |
| $f_{CH}$ | Exported flowrate to the CH from each plant |
| $f_D$ | The discharged flowrate of each plant |
| $f_F$ | Freshwater flowrate demanded by sinks |
| $f_{Imp}$ | Imported flowrate to sinks |
| $f_R$ | Reuse/recycle flowrate from sources to sinks |
| $f_W$ | Effluent flow rate from sources |
| $g_{CH}$ | The received flowrate of each plant from the CH |
| $l$ | Binary variable for export pipelines |
| $m$ | Binary variable for choosing the desired regenerator |
| $mrem$ | Total contaminant mass removed via regenerators |
| $Obj_{cost}$ | Objective function regarding the total annualized cost |
| $Obj_{Fresh}$ | Objective function regarding the freshwater consumption |
| $RC_t$ | The total annualized cost of the regenerator in the CH |
| $RR_t$ | The total removal ratio of the CH regenerator |
| $y$ | Binary variable for import pipelines |

*Input Parameters*

| | |
|---|---|
| $AWH$ | Annual working hours of the park |
| $AF$ | Annualized factor |
| $C_F$ | The concentration of contaminants in freshwater |
| $C_{SK}$ | The maximum concentration of contaminant in sinks |
| $C_{SR}$ | Maximum concentration of contaminant in sources |
| $D$ | Equivalent distance between participants and the CH |
| $EL$ | The economic life of the project |
| $F_{SK}$ | Limiting flowrate of sinks |
| $F_{SR}$ | Limiting flowrate of sources |
| $LB_{CH}$ | The lower bound of cross-plant flowrate |
| $P$ | The weight factor of different periods |
| $p$ | Cost parameter of pipelines |
| $q$ | Cost parameter of pipelines |
| $RR$ | Removal ratio of different regenerators |
| $RC$ | Regeneration unit cost of different regenerators |
| $u$ | The real interest rate during the lifespan of the project |
| $UB_{CH}$ | Upper bound of cross-plant flowrate |
| $\rho$ | Water density |
| $\upsilon$ | Water velocity in pipes |

# 1. Introduction

The quality and the value of sustainable development have recently become the focus of increased attention across the globe. In this regard, the determinants of sustainability within industrial zones are of great importance due to the impacts of industries on the quality of the environment, depletion of natural resources, and socio-economic issues (Susur et al., 2019, Taheri et al., 2021). As an example, industrial development is often linked to the use of high volumes of raw materials, energy, and freshwater (Huang et al., 2019). Achieving sustainability within the industrial sector calls for the more efficient industrial developments both in operation and design. These developments require powerful concepts since old methods fulfill the social, economic, and environmental needs, barely (Gómez et al., 2018). Consequently, there is a real necessity for different industries to ensure minimum natural resources consumption, while continuing high production intensities (Ramos et al., 2016, Taheri, et al., 2021).

The establishment of EIPs has been acknowledged as a significant step toward industrial sustainability because of EIPs' potential to reduce water and energy consumption, increase productivity and collaboration, and preserve the environment (Boix et al., 2015a). The EIP concept, which is directly linked to sustainable development (Kastner et al., 2015), seeks to join separate industries, geographically closed enough, in a cooperative manner so that exchanges of raw materials, by-products, energy, and utilities are maximized (Chertow, 2000). According to (Behera et al., 2012), a wide, accepted definition of EIP is "an industrial system of planned materials and energy exchanges that aims at minimizing energy and raw materials use, maximizing proficiency, minimizing waste, and creating sustainable economic, ecological and social relationships."

An efficient resource conservation approach is advantageous in the realm of the industrial sector as it boosts competitiveness through reduced operational costs and helps sustainable development. Process integration has been confirmed as a promising approach of maximizing potential resource conservation (Chew et al., 2008). Apart from the earliest enlargement in energy management approaches (Dhole and Linhoff, 1992; Shenoy, 1995; Smith and Ragan, 2005), much recent process integration, work on water network synthesis, has been stated in the literature, and there are several techniques for integrating water management (Boix et al., 2015b).

Water management is conducted by pinch analysis, mathematical optimization, or a combination of the two. The pinch analysis is an insight-based approach which targets the minimum freshwater and wastewater flowrates (Bandyopadhyay, 2006; El-Halwagi et al., 2003; Foo, 2009, 2007; Hallale, 2002; Li and Yao, 2004; Prakash and Shenoy, 2005). Although pinch analysis gives good results, it cannot be used in complex problems, for example, multi-contaminant systems and networks with complex operational constraints (Klemeš and Kravanja, 2013). To evade these shortcomings, mathematical optimization approaches are implemented. The mathematical optimization methods start by considering a water network superstructure followed by a desired objective function and technical/economic/environmental constraints. Then, water recovery and regeneration are achieved by simulating the desired network structure and the operational condition for a water network (Bagajewicz and Savelski, 2001; Hul et al., 2007; Karuppiah and Grossmann, 2008).

Apart from studying single water networks, a powerful tool to improve water recovery is through interplant water integration (IPWI), that is, the integration of water between different networks (Chew et al., 2010). Therefore, water-using processes could be grouped according to their geographical location or as separate plants. Accordingly, a water source has the potential to be reused/recycled to another network as a new source. IPWI problem first addressed in (Olesen and Polley, 1996). In this study, a pinch-based load table method has been developed for fixed load problems. Following the same approach, (Spriggs et al., 2004) used recovery pinch diagram for minimum flowrate targeting in the EIPs. The authors also proposed the use of a centralized mixing (termed as "centralized utility hub" in this current work) for conducting by-product or waste exchange among several process plants within an EIP. In 2008, (Chew et al., 2008) proposed two main structures for the design of EIPs: direct integration and indirect integration. The former distributes water directly through cross-plant pipelines, while the latter via a centralized utility. The study showed that the centralized scheme is not only more practical but also more flexible and cost beneficial compared to the direct water integration scheme.

Based on the proposed superstructures, many studies have been conducted. (Lovelady et al., 2009) developed an interception model and solved its associated nonlinear programming (NLP) problem. In their work, recycle, reuse and separation using interception devices are considered as possible approaches for managing wastewater. A cost objective function was given for the design of the water network structure with a minimum total annualized cost which included the costs of interception devices, freshwater, and wastewater treatment. (Boix et al., 2012) minimized the number of network connections, freshwater consumption, and regenerated water quantity via a multi-objective optimization model. (Rubio-Castro et al., 2012) considered the possibility of using inter-plant water in a park. The formulation was further extended into a property-based integration model in their later work (Rubio‐Castro et al., 2013). Their focus was on solving the mixed integer nonlinear programming (MINLP) to facilitate the computational burden of the EIPs water network design. (Liu et al., 2015) investigated the mixing possibilities across the plants to avoid redundant solutions by considering the limit of inter-plant stream pipes at the building stage of network superstructures. Both direct and indirect superstructures are explored in this study. The problem then mathematically formulated targeting the minimum freshwater consumption. Since the optimum for the entire EIP does not necessarily coincide with the solution satisfying the participants' goals, (Aviso et al., 2010) presented fuzzy mathematical programming to incorporate the benefits of networking related to the individual targets and objectives of the participants. (Lim and Park, 2009) developed a model to evaluate the water systems from two points of view: the life circle cost and the life circle assessment. (Alnouri et al., 2014) designed inter-plant water network taking into account pipe merging purpose.

Although EIPs have accepted sustainability objectives, the development of EIPs stays challenging due to a wide range of performance traits that need to be investigated in the system rather than focusing only on the economic, environmental and social aspects. The network performance may be vulnerable to the changes that take place within different periods, and the lack of flexibility makes EIPs suffer from the reluctance of different industries to participate (Zhang et al., 2017). However, most design methods that have been offered for the synthesis of EIPs are broadly based on a fundamental premise: the configuration of EIPs does not vary with time. Nevertheless, the EIP faces multiple changes in practice; new plants might be involved into the park, a company may undergo a stop in production, and the capacity of plants might also fluctuate for some reasons or purposes. The water network still has to be flexible to handle the new operating conditions.

To counteract the challenges above, the water network must be flexible. Some researchers have tried to model possible changes in EIPs and turned their attention to the multi-period problems in water network design. In this context, (Faria and Bagajewicz, 2011) considered a rise in the mass load of water-using units to address an increase in plant capacity regarding new water-using units planned to be added over time. (Burgara-Montero et al., 2013) presented a multi-period method to incorporate the seasonal alteration in the optimal treatment of industrial wastewaters. (Arredondo-Ramírez et al., 2015) optimized an agricultural water system with multi-period storage, reuse, and distribution. All of these studies reflected multi-period water integration merely within a plant, while for the problems regarding multiple plants in an EIP, (Liao et al., 2007) proposed a methodology for designing a multiple plant water network considering various operating scenarios, to improve the overall flexibility of the network. Following the same approach, (Aviso, 2014) examined different future situations which are expected to happen as the fate and plans of other participants are not entirely divulged. A robust optimization model is developed in this work to obtain the optimal network design, where the direct water allocation including minimum freshwater consumption was acquired subjecting to several likely scenarios within the EIP. In another study, (Kolluri et al., 2016) presented three robust optimization approaches to generate various solutions for resource conservation in EIPs. Their model takes various future events into account, though, the inside water use condition of each plant was not taken into consideration in their research. More recently, (Liu et al., 2017) proposed a mathematical formulation for considering predictable variations in the EIPs. The water network synthesis then results from predictable development and adjustment variations of the park.

Though attention has been paid to the challenge of multi-period water integration for EIPs, the flexibility of the EIPs' water networks is the concept that has not been given considerable attention as of yet- to the best of our knowledge. Our study aims at developing a water network through a multi-period planning. Based on the superstructure, a mathematical model is formulated for optimal design from both financial and environmental concerns. All of the researches mentioned above so far are only highlighted direct water integration while indirect superstructures are more practical and flexible. In this study, a centralized hub (CH) is employed to study indirect water reuse.

### 1.1. Problem Definition and Paper Contributions

In response to the gap discussed above, a generalized superstructure of an EIP is introduced, in this paper, and the mathematical optimization is then implemented to evaluate the optimal capacities and scheduling. On this basis, the main contributions of the paper and the novelty of this study can be summarized as follows:
- A refined superstructure is introduced for designing water networks of EIPs.
- A multi-period mathematical formulation is presented, and the final/optimum solution is carried out based on that.
- Trade-off between water network flexibility and performance, and vice versa, are discussed.

The rest of the paper is structured as follows. Section 2 presents the network superstructure and mathematical representation of the objective function and technical constraints based on the proposed superstructure. The model implementation is briefly discussed in section 3. Lastly, two cases are studied in part 4 followed by the conclusion in section 4.

## 2. The Proposed Superstructure and Associated Mathematical Formulation

### 2.1. Water Network with a Centralized Utility

The superstructure for indirect integration scheme is exhibited in Fig.1 in which a set of water networks of the fixed flow rates, with process sources and sinks, are given. Connected with pumps, splitters, distribution pipes, and mixers, these sinks and sources might be considered for water reuse/recycle as shown. Apart from being recycled/ reused to sinks in the local network (dotted line), to achieve further water recovery, process sources might also be sent to the CH (dashed line). The unused wastewater from the sources is then discharged to the environment (red line). After the existing process sources are completely utilized, freshwater is the external source to be used. The main advantage of using a CH is that it is more practical in handling a large number of water networks in the EIP scheme. A network containing a CH reduces the associated piping cost by gathering streams to be exported from each plant because geographical distances within a single water network are considerably shorter than between participants. Simultaneously, controllability difficulties are reduced by the reduction of unescapable fluctuations in stream flow rates and water concentrations, which, in turn, enhances water network flexibility and controllability (Chew et al., 2008; Ma et al., 2007). Besides, the purpose of discharging wastewater and exporting it to CH from the same location is to control the contaminant concentrations of water streams entering the CH, so as to avoid over-treatment of wastewater, along with the operation problems resulting from the variation of water quality.

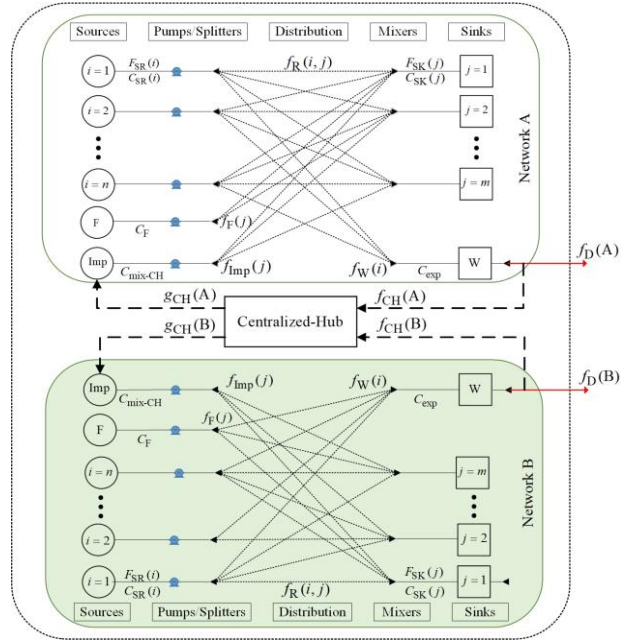

**Fig.1.** Superstructure for water network of EIPs via CH

### 2.2. Mathematical Formulation

The mathematical model for the proposed superstructure is expressed as follows. Typically, the water network model involves the objective functions, mass balance and flowrate constraints over each sink, source, and regenerator, and logical constraints demonstrating stream use specifications in different periods. The multi-period model is extended from the one-period design formulas suggested in the study of (Chew et al., 2008). To achieve the multi-period model, the dimension $t$ is introduced into all equations to further extend the formulas to different periods. Mass balances are involved in stating direct water use respecting in-plant and indirect water reuse regarding inter-plant regeneration operation. To indicate the usage of cross-period pipeline connections, rational constraints are added.

To characterize each stream and operation in the model, $i, j, k, t$ representing water source, water sink, plant, and time, respectively, are used in formulas. Under the premise of ignoring the mass losses in delivery and operation, the flow rate balance over the source unit $(i,t)$ is presented in eq (1). The entire flowrate of each source $(F_{SR})$ equals the sum of that allocated to sinks, where variable $(f_R(i,j,t))$ represented reuse/recycle flowrate from source $i$ to sink $j$ in the period $t$, and the discharge/export sink $(f_W(i,t))$.

$$F_{SR}(i,t) = \sum_{j=1}^{m} f_R(i,j,t) + f_W(i,t) \tag{1}$$

$$\forall i \in \{1,2,...,n\}, \forall t \in \{1,2,...,r\}$$

The received water of each sink unit $(F_{SK}(j,t))$ comes from sources, freshwater source $(f_F(j,t))$, and the imported water from CH $(f_{Imp}(j,t))$. Equation (2) describes the flow rate balance around the mixing points, and eq (3) states the maximum acceptable contaminant on the inlet in each unit.

$$F_{SK}(j,t) = \sum_{i=1}^{n} f_R(i,j,t) + f_F(j,t) + f_{Imp}(j,t) \tag{2}$$

$$\forall j \in \{1,2,...,m\}, \forall t \in \{1,2,...,r\}$$

$$\sum_{i=1}^{n} \left(f_R(i,j,t) \times C_{SR}(c,i,t)\right) + \left(f_F(j,t) \times C_F(c)\right)$$
$$+ \left(f_{Imp}(j,t) \times C_{\text{mix-CH}}(c,t)\right) \leq F_{SK}(j,t) \times C_{SK}(c,j,t) \tag{3}$$

$$\forall c, \forall j \in \{1,2,...,m\}, \forall t \in \{1,2,...,r\}$$

In this equation, $C_{SR}(c,i,t)$ is the sources' contamination level, $C_F(c)$ is the freshwater contamination which is set to zero according to (Chew et al., 2008), $C_{mix\text{-}CH}(c,t)$ is the contamination of the outlet streams from the CH, and $C_{SK}(c,j,t)$ is the maximum acceptable contaminant in sinks.

Equation (4) exhibits the flowrate balance of the splitter which imports water from the CH. In any given period, the sum of imported water of sinks in one plant equals to the flowrate sent by the CH to that participant $(g_{CH}(k,t))$.

$$g_{CH}(k,t) = \sum_{j=1}^{m} f_{Imp}(j,t) \quad (4)$$
$$\forall k \in \{1,2,...,K\}, \forall t \in \{1,2,...,r\}$$

Equation (5) states that the exported water from sources is split into two streams after mixing: be discharged $(f_D(k,t))$ or sent to the CH $(f_{CH}(k,t))$.

$$\sum_{i=1}^{n} f_W(i,t) = f_{CH}(k,t) + f_D(k,t) \quad (5)$$
$$\forall k \in \{1,2,...,K\}, \forall t \in \{1,2,...,r\}$$

With binary variables $l(k,t)$ and $y(k,t)$ signifying the present of cross-plant pipelines, eqs (6) and (7) set, respectively, the upper and lower bounds of the cross-plant flowrates to and from the CH. When $l(k,t)$ or $y(k,t) = 1$ the pipelines exist on a particular period. $LB_{CH}$ and $UB_{CH}$ are also indicating the lower and upper quantity of feasible water flow rates in the cross-plant pipes.

$$LB_{CH} \times l(k,t) \leq g_{CH}(k,t) \leq UB_{CH} \times l(k,t) \quad (6)$$
$$\forall k \in \{1,2,...,K\}, \forall t \in \{1,2,...,r\}$$

$$LB_{CH} \times y(k,t) \leq f_{CH}(k,t) \leq UB_{CH} \times y(k,t) \quad (7)$$
$$\forall k \in \{1,2,...,K\}, \forall t \in \{1,2,...,r\}$$

Also, the level of contamination in the discharged and exported streams can be formulated as below:

$$\sum_{i=1}^{n} f_W(i,t) \times C_{SR}(c,i,t) = \left(f_{CH}(k,t) \times C_{exp}(c,k,t)\right) + \left(f_D(k,t) \times C_{exp}(c,k,t)\right) \quad (8)$$
$$\forall c, \forall k \in \{1,2,...,K\}, \forall t \in \{1,2,...,r\}$$

where $C_{exp}(c,k,t)$ represents the level of contamination in discharged/exported streams of each plant in any given period. Under the assumption that the CH only change the level of contaminants, eq (9) formulates the flowrate balance in the centralized utility.

$$\sum_{k=1}^{K} f_{CH}(k,t) = \sum_{k=1}^{K} g_{CH}(k,t), \forall t \in \{1,2,...,r\} \quad (9)$$

Despite the flowrates, contaminants are highly affected by the CH. Equation (10) formulates this change as:

$$\sum_{k=1}^{K} f_{CH}(k,t) \times C_{exp}(k,t) = \sum_{k=1}^{K} g_{CH}(k,t) \times C_{mix\text{-}CH}(t) + mrem(c,t) \quad \forall t \in \{1,2,...,r\} \quad (10)$$

where $mrem(c,t)$ is the mass removed from the flowrates in the centralized utility. The CH entails regenerators to reduce the contamination level of exported streams. The quantity of removed mass from wastewaters is highly correlated with the investment costs of the operating regenerator; the higher a regenerator unit can remove contaminants, the more it costs. To characterize each regenerator and its operation in the model, $s$ representing various regenerators is used in formulas. In this regard, the regenerator units are designated with their removal ratio coefficients $(RR)$ and their associated costs $(RC)$. Table 1 describes the characteristics of several regenerators according to (Liu et al., 2017).

**Table 1**
Cost for wastewater removing for different regenerators

| S | RR | $RC\ [\$/kg\ \text{removed}]$ |
|---|---|---|
| 1 | 0.1 | 0.540 |
| 2 | 0.2 | 0.695 |
| 3 | 0.3 | 0.850 |
| 4 | 0.4 | 1.005 |
| 5 | 0.5 | 1.160 |
| 6 | 0.6 | 1.460 |
| 7 | 0.7 | 1.760 |
| 8 | 0.8 | 2.060 |

The total mass removed from wastewater can be written as:

$$mrem(t) = \left(\sum_{k=1}^{K} f_{\text{CH}}(k,t) \times C_{\exp}(c,k,t)\right) \times RR_{\text{t}}(c,t) \quad (11)$$

$$\forall c, \forall t \in \{1, 2, ..., r\}$$

where $RR_{\text{t}}(c,t)$ is the total removal ratio which is calculated via eq (12).

$$RR_{\text{t}}(c,t) = \sum_{s=1}^{S} m(s,t) \times RR(c,s), \forall c, \forall t \in \{1, 2, ..., r\} \quad (12)$$

where $m(s,t)$ is the binary variable indicating the level of regeneration in operation. Since only one regenerator will operate in the CH, we have:

$$\sum_{s=1}^{S} m(s,t) = 1, \forall t \in \{1, 2, ..., r\} \quad (13)$$

To keep continuous operation for regenerators and cross-plant pipelines, the unit or pipes should be in use since its first setting; otherwise, its service decided by the contemporary design for a particular period. Thus, logical constraints in eqs (14) to (16) are introduced, respectively, for the regenerators within the CH, cross-plant pipelines from CH to participants and plants to CH. These constraints call for

$$m(s,t+1) - m(s,t) \geq m(s,t) - 1, \forall t \in \{1, 2, ..., r\} \quad (14)$$
$$y(k,t+1) - y(k,t) \geq y(k,t) - 1, \forall t \in \{1, 2, ..., r\} \quad (15)$$
$$l(k,t+1) - l(k,t) \geq l(k,t) - 1, \forall t \in \{1, 2, ..., r\} \quad (16)$$

For instance, $l(k,t+1)$ would be equal to 1 if the pipeline existed in the prior period, and it could be equal to 0 or 1 otherwise.

The minimum weighted total annualized cost is objective represented by eq (17). Where, $P(t)$ is the weight factor that indicates the influence of period $t$ to the objective.

$$\text{Obj}_{\text{cost}} = \sum_{t=1}^{r} P(t) \times \left(C_{\text{inv}}(t) + C_{\text{op}}(t)\right) \quad (17)$$

where the right side of eq (17) is divided into two terms: $C_{\text{inv}}(t)$ is the investment cost for pipelines and $C_{\text{op}}(t)$ is the operating cost including regeneration operation, freshwater consumption, and wastewater treatment in any given period. The investment and operation costs are detailed in eqs (18) and (19) according to (Chew et al., 2008; Liu et al., 2017).

$$C_{inv}(t) = D \times AF \times p \times \left(\sum_{t=1}^{r}\sum_{k=1}^{K}\left(\frac{f_{\text{CH}}(k,t) + g_{\text{CH}}(k,t)}{3600\rho v} + \frac{q}{p}y(k,t) + \frac{q}{p}l(k,t)\right)\right) \quad (18)$$

$$C_{op}(t) = AWH \times 10^{-6} \times \left(\begin{array}{l}\sum_{j=1}^{m}\left(f_{\text{F}}(j,t) \times \pi_{\text{F}}\right) + \\ \left(RC_{\text{t}}(t)\right) + \sum_{k=1}^{K}\left(f_{\text{D}}(K,t) \times \pi_{\text{D}}\right)\end{array}\right) \quad (19)$$

starting with eq (18), $D$ is an equivalent distance between the CH and plants, assumed 100m throughout this study. $AF$ is the annualized factor given in eq (20), $p$ and $q$ are cost parameters for pipelines which are respectively equal to 7200 and 250, $\rho$ is water density equals to 1000 kg.m$^{-3}$, and $v$ is the flowrate velocity equals to 1 m.s$^{-1}$. $AWH$ in eq (19) represents annual working hours of the EIP which is assumed 8000h during this study.

$$AF = \frac{u(1+u)^{EL}}{(1+u)^{EL} - 1} \quad (20)$$

where $u$ is the fractional interest rate per year, and $EL$ is the economic life of the project. Finally, the associated total annualized cost of the regenerators can be written as:

$$RC_t(t) = \sum_{s=1}^{S} \left( m(s,t) \times RC(s) \times mrem(t) \times 10^{-3} \right) \quad (21)$$

$$\forall t \in \{1, 2, ..., r\}$$

## 3. Implementation

To solve the proposed MINLP model, GAMS is used, and the model is solved with its built-in solvers. To be clear, the achieved results might not be globally optimum due to the complexity of the established MINLP model. However, a reformulation-linearization approach is implemented according to the work of (Quesada and Grossmann, 1995) to alleviate the complexity of the problem. It is worthy of note that the obtained solutions and conclusions would not be affected if the answers were just locally optimal.

## 4. Case Studies

In this segment, two EIPs from (Chew et al., 2008) and (Rubio-Castro et al., 2012) are considered in two scenarios to illustrate the results. Case 1 is the base case in which the EIPs faces no changes in their configuration during the lifespan. This scenario is presented in order to demonstrate the effectiveness of the proposed superstructure and act as the base case for further comparisons. Networks of this case are optimized aiming at the minimum annualized cost expectation. Case 2 is a three-period problem concerning new plant construction scenarios in the EIP. The information of the EIP's participants including the concentrations of sources and their offered flowrates and the required flowrates of sinks together with their concentration limits are given in Table 2 and Table 3.

**Table 2**
The participants' data of EIP-1

| water network | source SR($i$) | flowrate $F_{SR}$ (t/h) | concentration $C_{SR}$ (ppm) | sink SK($j$) | flowrate $F_{SK}$ (t/h) | concentration $C_{SK}$ (ppm) |
|---|---|---|---|---|---|---|
| A | 1 | 20.00 | 100 | 1 | 20.00 | 0 |
|  | 2 | 66.67 | 80 | 2 | 66.67 | 50 |
|  | 3 | 100.00 | 100 | 3 | 100.00 | 50 |
|  | 4 | 41.67 | 800 | 4 | 41.67 | 80 |
|  | 5 | 10.00 | 800 | 5 | 10.00 | 400 |
| B | 6 | 20.00 | 100 | 6 | 20.00 | 0 |
|  | 7 | 66.67 | 80 | 7 | 66.67 | 50 |
|  | 8 | 15.63 | 400 | 8 | 15.63 | 80 |
|  | 9 | 42.86 | 800 | 9 | 42.86 | 100 |
|  | 10 | 6,67 | 1000 | 10 | 6,67 | 400 |
| C | 11 | 20.00 | 100 | 11 | 20.00 | 0 |
|  | 12 | 80.00 | 50 | 12 | 80.00 | 25 |
|  | 13 | 50.00 | 125 | 13 | 50.00 | 25 |
|  | 14 | 40.00 | 800 | 14 | 40.00 | 50 |
|  | 15 | 300.00 | 125 | 15 | 300.00 | 100 |

**Table 3**
The participants' data of EIP-2

| water network | source SR($i$) | flowrate $F_{SR}$ (t/h) | concentration $C_{SR}$ (ppm) | sink SK($j$) | flowrate $F_{SK}$ (t/h) | concentration $C_{SK}$ (ppm) |
|---|---|---|---|---|---|---|
| A | 1 | 450.00 | 15 | 1 | 450 | 11 |
|  | 2 | 320.00 | 17 | 2 | 320 | 13 |
|  | 3 | 300.00 | 50 | 3 | 300.00 | 37 |
| B | 4 | 200.00 | 20 | 4 | 200.00 | 17 |
|  | 5 | 150.00 | 22 | 5 | 150.00 | 16 |
|  | 6 | 75.00 | 35 | 6 | 75.00 | 11 |
| C | 7 | 120.00 | 48 | 7 | 120.00 | 25 |
|  | 8 | 600.00 | 49 | 8 | 600.00 | 20 |
|  | 9 | 300.00 | 50 | 15 | 300.00 | 13 |

## 4.1. Case 1: One-period Case with All participants, Targeting the Minimum Annualized Cost

In this case, the MINLP model is solved using the objective function in eq (17) to minimize the total annualized cost, subject to constraints in eqs (1) to (16). Since there is only one period, the weight factor ($P(t)$) is equal to 1. This case is formulated into a model consisting of 159 continuous variables, 14 binary variables, and 169 constraints for EIP-1 and 98 continuous variables, 6 binary variables, and 88 constraints for the EIP-2. Fig. 2 shows the optimal developing scheme of the networks for the minimum total annualized cost. The overall annualized cost of the EIP-1 network is calculated 1.24 million which is 4.17% less than what was obtained in the study of (Chew et al., 2008). Also, the objective function is reduced from 0.66 to 0.62 m$/year for EIP-2. The total wastewater from the networks (which are 142.90 t/h for EIP-1 and 66.14 t/h for EIP-2) shows a subtle reduction (-1.92% and -24.13%) for EIP-1 and EIP-2, respectively. In addition, the freshwater flowrate is decreased from the 145.72 t/h in the study of (Chew et al., 2008) to 142.89 t/h in the present work for EIP-1.

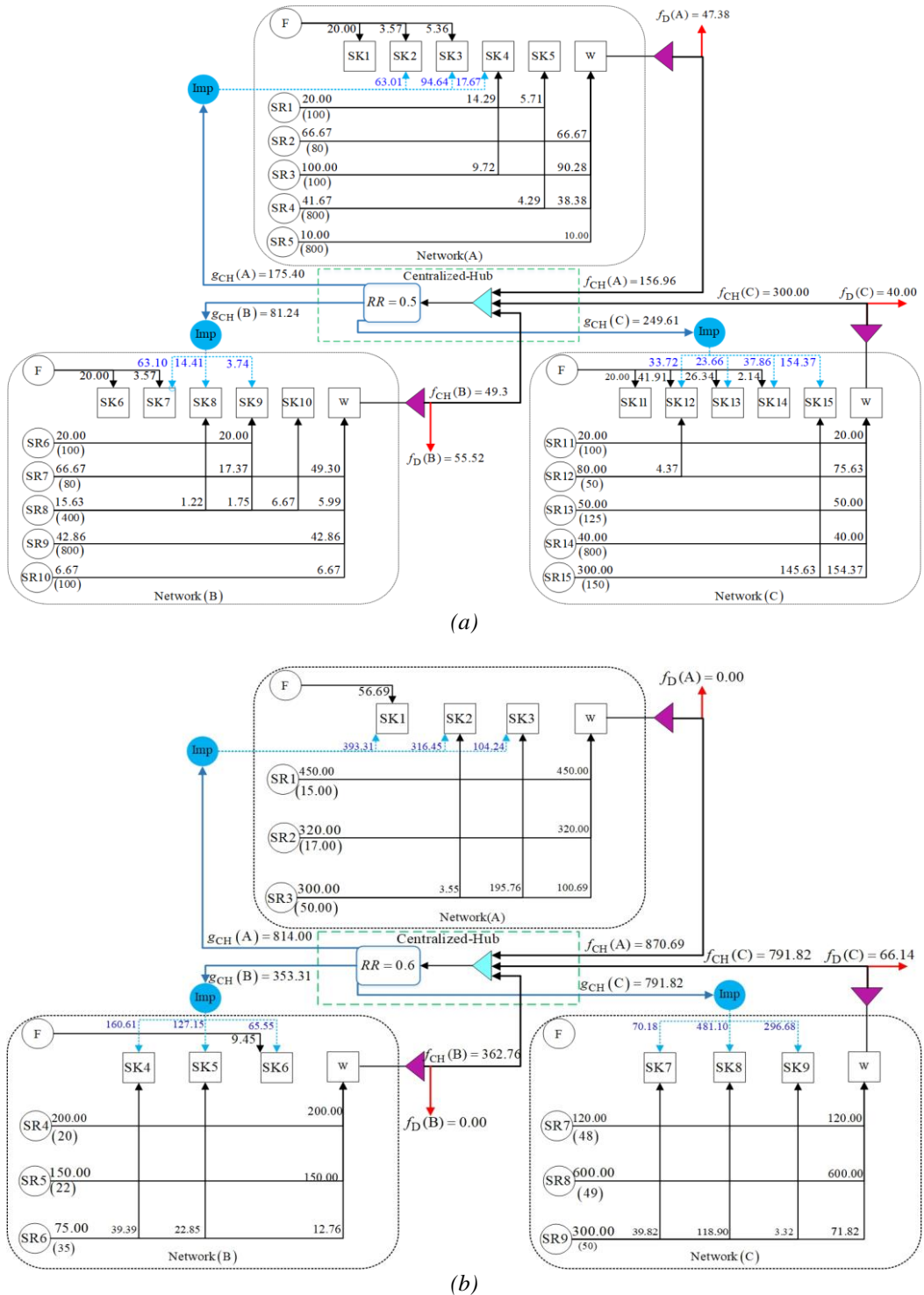

**Fig.2.** One-period design of the water networks (optimized for minimum total annualized cost)
(a) EIP-1 (b) EIP-2

*4.2. Case 2: Three-period Design with New Participant's entry, targeting the minimum annualized cost*

Case 2 seeks the most economical design for the networks through trading off amongst the costs for freshwater, wastewater treatment and regeneration operation, and pipes investment. In this case, it is assumed that plant A is the participant and plants B and C enter into the park in the 3rd and 5rd year of the EIP's interesting horizon, respectively. Some assumptions are considered to simplify the synthesis: (1) any change in network's configuration is accomplished at the start of the corresponding period, instantly, (2) the weight factor of different periods are equal. The optimal working condition of the EIP-1 during different periods are shown in Fig. 3 concerning the economic objective function. A regenerator with the removal ratio of ($RR$=0.4) is applied, and six pipelines are in operation to send streams from the CH to plants A and B and vice versa. The weighted total annualized cost is equal to 1.39 million.

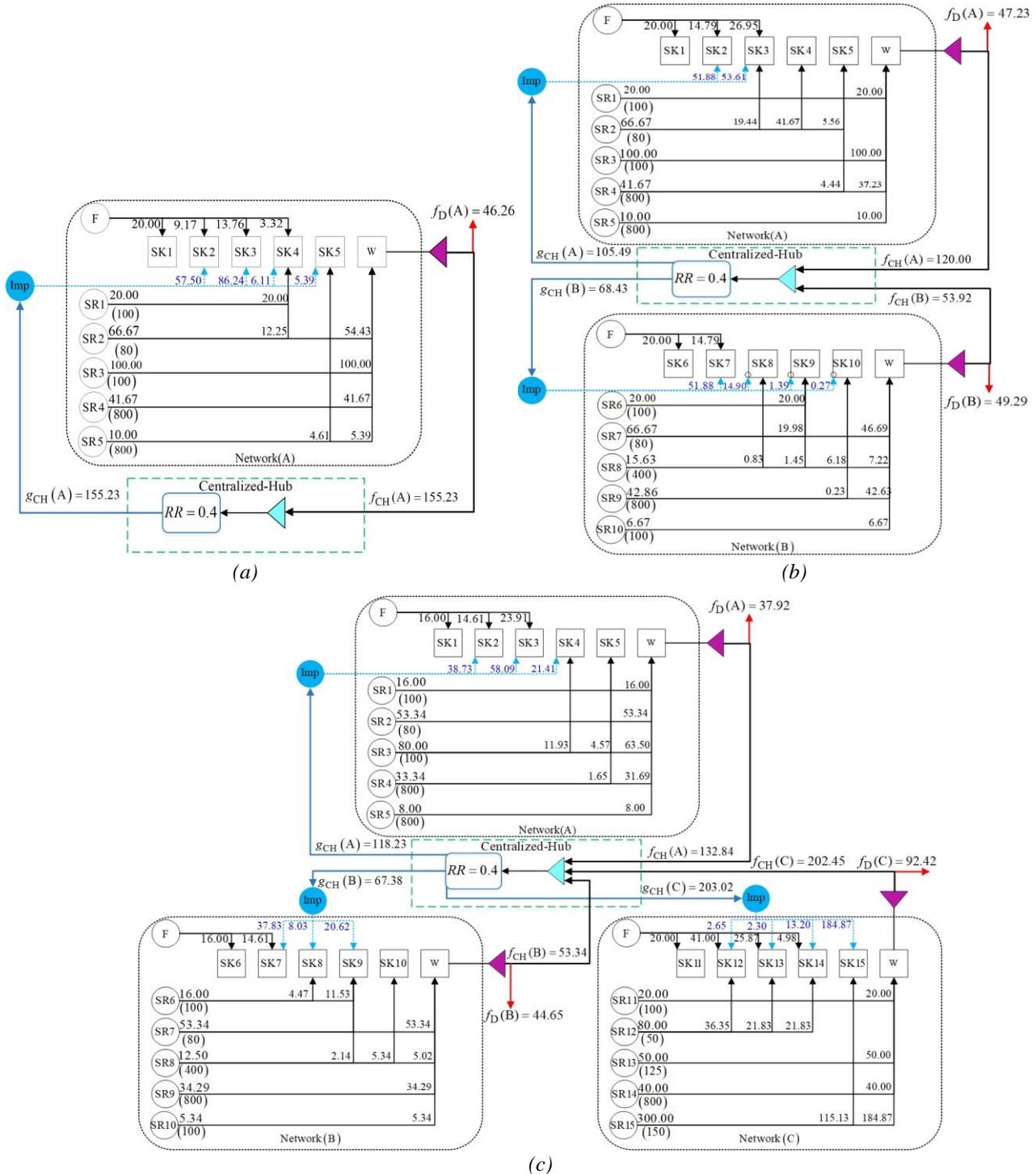

**Fig.3.** Three-period design of the water network (optimized for minimum total annualized cost). EIP-1
*(a) t1  (b) t2  (c) t3*

Upon entering plant C in the EIP, some changes are made in the park's arrangement and flows. The associated model contains 290 continuous variables, 28 binary variables, and 295 constraints. The participants' networks and the corresponding CH pipelines linking with plants are determined after the overall design over two periods is specified. Nevertheless, the EIP's connections may not be constantly in operation due to the different water use specifications in different periods. As more participants are involved, the degree of pipes' utilization rises and vacancy rate declines. This demonstrates that the model belongs to a multi-period problem rather than a single-period problem. Fig. 3 (c) exhibits the optimized scheme of the network in this multi-period case. To operate continuously, the existing pipes between the CH and plants and the regenerator with ($RR$=0.4) should be in operation because of their setting in the period $t1$. As indicated in this scheme, the total freshwater flowrate is 155.28 t/h which shows a significant increase compared to case 1. Also, the total wastewater flowrate is equal to 175.00 t/h. The total annualized cost is increased from 1.24 million in case 1 to 1.39 million in this case, considering equal weight factor for periods $t1$, $t2$ and $t3$. Following the same approach, the minimum weighted freshwater consumption is equal to 155.28 t/h.

The optimal working condition of the EIP-2 for different periods are presented in Fig. 4. A regenerator with the removal ratio of ($RR$=0.6) is used. The weighted total annualized cost is equal to 1.08 million.

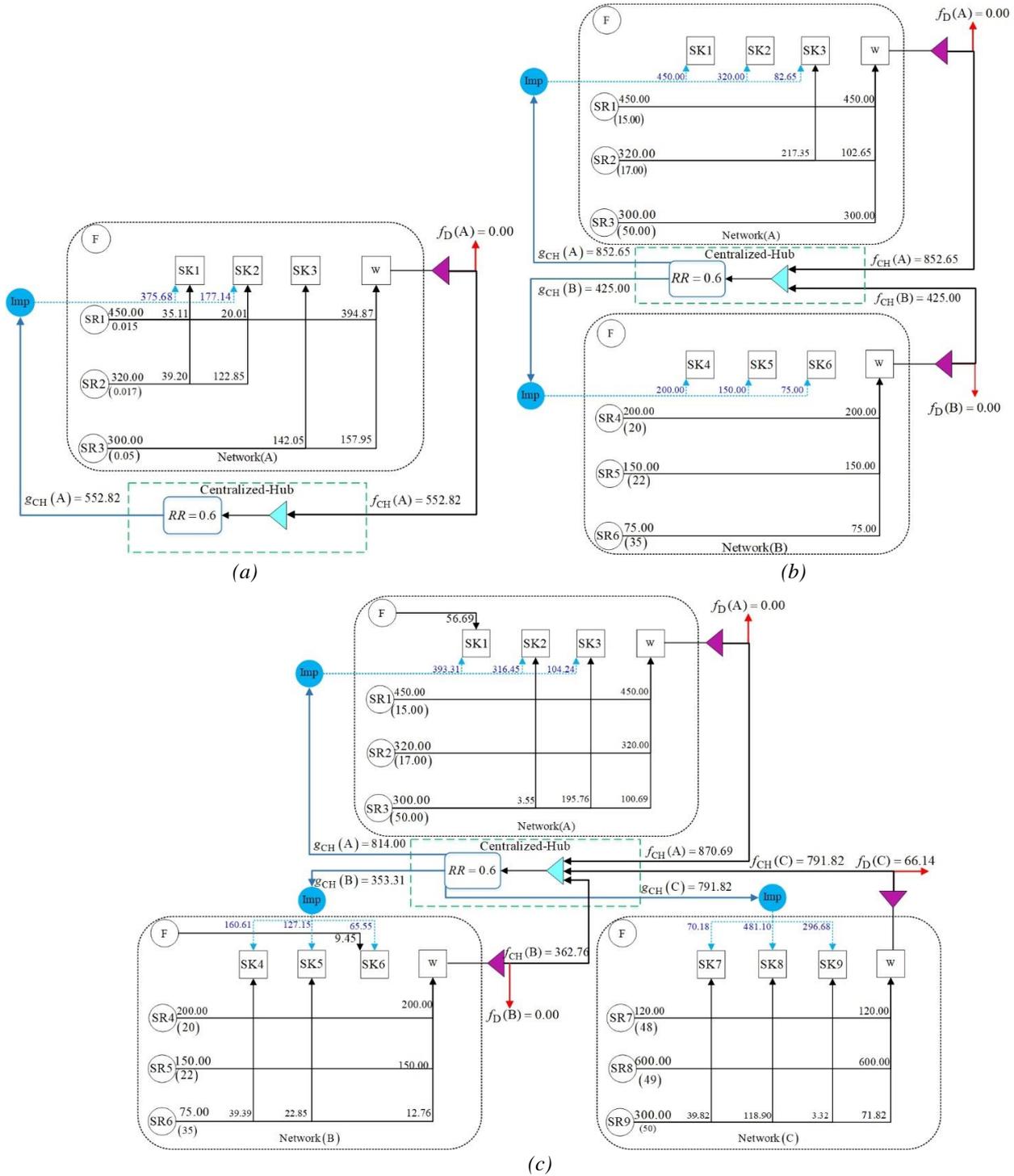

**Fig. 4.** Three-period design of the water network (optimized for minimum total annualized cost), EIP-2

Table 3 summarizes the results obtained in the proposed cases including the total freshwater and wastewater flow rates, the removal ratio of the implemented regenerator, the minimum weighted expected annualized cost, and the overall NCI of the network. As the table indicated, the gradual emerging of plants has increased the total freshwater consumption and costs in the EIP.

**Table 3**
Results for water network design of different case studies for two EIPs

| Model | freshwater flowrate (t/h) | wastewater flowrate (t/h) | removal ratio of the CH | weighted annualized cost (million $/year) | freshwater cost (million $/year) | wastewater cost (million $/year) | cross-plant pipeline capital cost(million $/year) |
|---|---|---|---|---|---|---|---|
| Case 1-EIP-1 | 142.89 | 142.90 | 0.5 | 1.24 | 0.58 | 0.58 | 0.02 |
| Case 2-EIP-1 | 155.28 | 175.00 | 0.4 | 1.39 | 0.63 | 0.71 | 0.02 |
| Case 1-EIP-2 | 56.69 | 66.14 | 0.6 | 0.62 | 0.23 | 0.28 | 0.02 |
| Case 2-EIP-2 | 88.73 | 150.35 | 0.6 | 1.08 | 0.36 | 0.61 | 0.02 |

## 5. Conclusion

An interplant water integration superstructure with a centralized utility together with a regeneration unit is introduced in this study with its purpose to improve water recovery through the EIPs. Furthermore, a multi-period model for the water network is developed and solved aiming at minimizing the weighed freshwater flowrate and weighed the total annualized cost. The proposed superstructure operates with lower freshwater flowrate and overall cost compared to the networks presented in the previous studies. Considering various operation scenarios will make the best planning of the water network for the EIP. Besides, a resilience indicator, adopted for EIPs, used and the overall resilience of the park are evaluated for four different scenarios. Generally, the resilience indicator is decreased with the involvement of new participants. The dynamic development of water networks given in figures has confirmed that to motivate the participation of companies, further studies are required to evaluate different aspects of the water network designs for EIPs, apart from considering only the costs and water consumption. Also, the resilience indicator can be involved in multi-objective models to achieve more resilient networks.